\documentclass[12 pt]{article}
\usepackage[utf8]{inputenc}
\usepackage[myheadings]{fullpage}
\DeclareUnicodeCharacter{0301}{\hspace{-1ex}\'{ }}
\usepackage{fancyhdr}
\usepackage{lastpage}

\usepackage{graphicx, wrapfig, subcaption, setspace, booktabs}
\usepackage[T1]{fontenc}

\usepackage[font=small, labelfont=bf]{caption}
\usepackage{fourier}
\usepackage[protrusion=true, expansion=true]{microtype}

\usepackage{amsmath,amssymb}
\usepackage{float}
\usepackage{graphicx}
\usepackage{wrapfig}
\usepackage[colorinlistoftodos]{todonotes}
\usepackage[colorlinks=true, allcolors=black]{hyperref}

\usepackage[english]{babel}
\usepackage{csquotes}

\usepackage[a4paper,top=2cm,bottom=1.5cm,left=2cm,right=2cm,marginparwidth=2.5cm]{geometry}

\newtheorem{Theorem}{Theorem}[section]
\newtheorem{Remark}{Remark}[section]

\newtheorem{Lemma}{Lemma}[section]

\newtheorem{Definition}{Definition}[section]
\newtheorem{Example}{Example}
\numberwithin{equation}{section}
  
\begin{document}
\begin{center}
{\textbf {\Large On astheno-K\"ahler manifolds}}\\[2pt]
{\bf Punam Gupta}\footnote{
School of Mathematics, Devi Ahilya Vishwavidyalaya, Indore-452 001, M.P. 
INDIA\newline
Email: punam2101@gmail.com} and {\bf Nidhi Yadav}\footnote{
Department of Mathematics \& Statistics, Dr. Harisingh Gour Vishwavidyalaya,
Sagar-470 003, M.P. INDIA \newline
Email: nidhiyadav.bina@gmail.com} 
\end{center}

\noindent
{\bf Abstract}
This survey explores a range of classical findings and recent developments related to our understanding of astheno-K\"ahler manifolds. Furthermore, we provide various examples of astheno-K\"ahler manifolds and analyze the challenges associated with their existence.

\noindent {\bf MSC: } 53C55, 32C35.
\vskip.2cm
  \noindent
  Keywords: Complex manifolds, astheno-K\"ahler metrics,  cohomology, strong KT metric,  Gauduchon metrics.

\section{Introduction}

During the last $30$ years, astheno-K\"ahler manifolds have attracted significant interest among geometers worldwide, showing potential relevance in both differential and algebraic contexts. The first four sections present a brief overview of recent advances in the study of astheno-K\"ahler manifolds, accompanied by examples and counterexamples. The final section delves into specific details regarding the estimates of equations and their solutions within the astheno-K\"ahler manifolds framework. 
\section{Preliminaries}
An almost complex structure on a real differentiable manifold $M$ of complex dimension $m$
 is a tensor field $J$ which at each point $x \in M$ is an endomorphism of the tangent space $T_{x}M$ such that $J^{2} = -I$, where $I$ denotes the identity transformation of
 $T_{x}M$.  A manifold with such a structure is called an almost complex manifold. An almost complex structure is called a complex structure if and only if $J$ has no torsion \cite{Yano}. 

 \begin{Definition}
  {\rm \cite{nijen}} Let $(M,J)$ be an almost complex manifold, then a vector-valued $2$-form $N(J, J)$ associated with $J$ is given by 

$$
N(J, J)(X, Y)=-[J X, J Y]+J[X, J Y]+J[JX, Y]+[X, Y],
$$
where $X$ and $Y$ are vector fields. This form is called the torsion tensor, or the Nijenhuis tensor of the almost-complex structure $J$. $J$ is said to be integrable if $N = 0$, then $J$ is a complex structure and $(M, J)$ is a complex manifold.

\end{Definition}
 
 An almost complex manifold $M$ \cite{scyau} together with a compatible
Riemannian metric $g$, that is, $g(JX, JY ) = g( X, Y)$
for all vector fields $X$, $Y$ on $M$, is called an almost Hermitian manifold and metric $g$ is called almost Hermitian metric. A complex manifold $M$ \cite{scyau} together with a compatible
Riemannian metric $g$ is called a Hermitian manifold and metric $g$ is called a Hermitian metric. 

The alternating $2$-form
 $$\Omega(X, Y ) = g(JX, Y)$$
is called the associated K\"{a}hler form. We can retrieve $g$ from $\Omega$
$$g( X, Y ) = \Omega(X, JY ).$$
\noindent If $\Omega $ is closed, then $(M,g)$ is known as a K\"{a}hler
manifold  and $g$ is a K\"{a}hler metric.

A quasi-K\"ahler structure is an almost Hermitian structure whose K\"ahler form $\Omega$ satisfies $$(d\Omega)^{(1,2)} = \bar{\partial} \Omega =0.$$

The term astheno-K\"ahler manifold (astheno is a Greek word, which means weak) was first coined in 1993 by Jost and Yau in the remarkable paper \cite{jsty,jst}.

\begin{Definition} {\rm\cite{jsty,jst}} A Hermitian manifold $(M,J,g)$ of complex dimension $m$ is called an astheno-K\"ahler manifold if it carries a fundamental $2$-form (K\"ahler form) $\Omega$ satisfying

 $$\partial \bar{\partial} \Omega^{m-2}=0,$$ where $\partial $
 and $\bar{\partial}$ are complex exterior differentials and $\Omega^{m-2}=\underset{m-2\quad \mathrm{times}}{\Omega \wedge \ldots\wedge \Omega}$.
\end{Definition}
The above condition is automatically satisfied for $m=2$. Thus, any Hermitian metric on a complex surface is astheno-K\"ahler.

\begin{figure}[h!]
\begin{center}
\includegraphics[width=.8\textwidth]{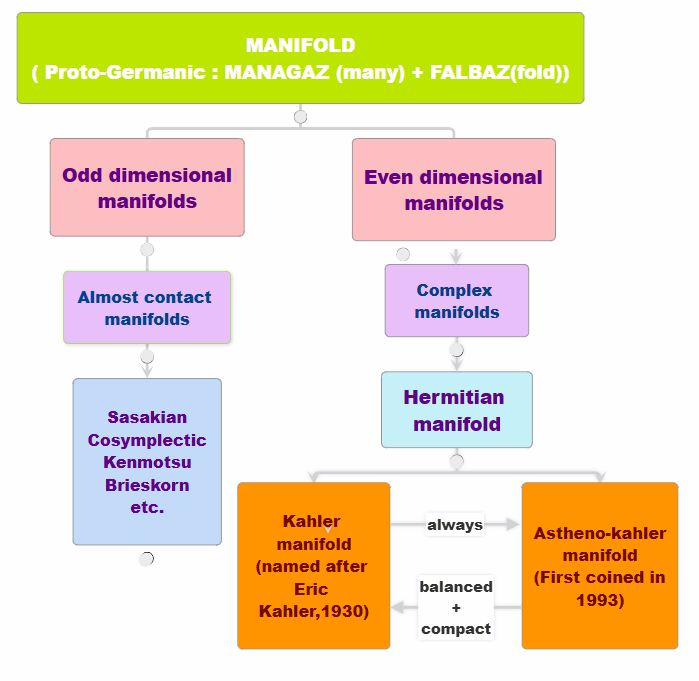}
\caption{flow chart}
\end{center}
\end{figure}

\begin{Definition}
   {\rm \cite{km}} Let $\nabla$ be a connection on a smooth manifold $M$. The torsion tensor $T$ is a $(1,2)$-tensor (alternatively, a vector-valued 2-form) defined as:

$$
T(X, Y)=\nabla_X Y-\nabla_Y X-[X, Y]
$$
where, $X, Y$ are vector fields,$\nabla_X Y$ denotes the covariant derivative of $Y$ along $X$ and $[X, Y]$ is the Lie bracket of $X$ and $Y$.\\
\noindent
For a connection $\nabla$ with torsion $T$, a torsion $1$-form $\tau$ can be defined in terms of certain contraction:

$$
\tau(X)=\operatorname{trace}(Y \mapsto T(X, Y))
$$

\noindent
This contraction gives a scalar-valued measure of the torsion in the direction of $X$. \\
The torsion $1$-form is related to the torsion of a connection on a manifold, which quantifies how the connection deviates from being flat. It characterizes how vectors are transported with respect to that connection. 
\end{Definition}




\begin{Definition}
{\rm \cite{Streets}} Let $M$ be an $m$-dimensional Hermitian manifold. The Hermitian structure $(J,g)$ is called strong K\"ahler with torsion (strong KT or SKT)  and $g$ is called strong KT metric or SKT metric if $ \partial\Bar{\partial}\Omega = 0.$
\end{Definition}
For $m=3$, astheno-K\"ahler structure means a strong KT metric.

\begin{Definition}
    Let a Hermitian manifold $(M,J,g)$ of complex dimension $n$ carry a fundamental $2$-form (K\"ahler form) $\Omega$ defined as $\Omega(X, Y ) = g(JX, Y)$. The Lee form, $\theta$, is the unique $1$-form defined by:

$$
\theta=-\frac{1}{n-1} \delta \Omega,
$$where $\delta$ is the codifferential operator (dual to the exterior derivative $d$ ).
In this context, the Lee form quantifies the extent to which the almost Hermitian structure deviates from being Hermitian or K\"ahler.

\end{Definition}

\begin{Definition}
   {\rm \cite{kobayashi}} Let $(M, g, J)$ be a Hermitian manifold. The Hermitian connection is a unique connection $\nabla$ on $M$ that satisfies $\nabla g=0$ (metric compatibility) and $\nabla J=0$ (compatibility with the complex structure). The curvature of the Hermitian connection is encoded in the Hermitian curvature tensor $R^{\nabla}$. A Hermitian manifold is said to be Hermitian flat if $
R^{\nabla}=0,
$
where $$R^{\nabla}(X, Y) Z=\nabla_X \nabla_Y Z-\nabla_Y \nabla_X Z-\nabla_{[X, Y]} Z.$$
\end{Definition}

 \begin{Definition}
{\rm \cite{gaud}} Let $M$ be an $m$-dimensional Hermitian manifold, then the Hermitian structure $(J,g)$ is called standard or Gauduchon and $g$ is called Gauduchon metric if $ \partial\Bar{\partial}\Omega^{m-1} = 0$. Equivalently, if the Lee form of $\Omega^{m-1}$ is co-closed. 
\end{Definition}
Fino et al. {\rm \cite{at}} mentioned if $m = 4$, a Hermitian metric that is simultaneously strong KT and astheno-K\"ahler metric, it must also be Gauduchon.

\vspace{0.5cm}
 \noindent
Fu et al. \cite{fu} introduced and investigated the generalization of Gauduchon metrics called the $k$-th Gauduchon. 
 \begin{Definition}
Let $M$ be an $m$-dimensional Hermitian manifold and $k$ be an integer such that $1\leq k \leq m-1$, then the Hermitian structure $(J, g )$ is called $k$-th Gauduchon and $g$ is called $k$-th Gauduchon metric if $$\partial \bar{\partial} \Omega^{k}\wedge \Omega^{m-k-1}=0.$$
 \end{Definition}

\begin{Definition}
  {\rm \cite{km}} A Hermitian manifold is said to be balanced if the torsion $1$-form of its Hermitian connection vanishes everywhere or $d \Omega^{m-1}=0$. Equivalently, the Hermitian structure is said to be balanced if its Lee form vanishes.
\end{Definition}

\noindent Now, it is natural to ask, when astheno-K\"ahler manifold becomes K\"ahler manifold?

\noindent Matsuo and Takahashi \cite{km} gave a significant result and showed that every compact balanced astheno-K\"ahler manifold is K\"ahler. In particular, compact Hermitian-flat astheno-K\"ahler manifolds are K\"ahler. 

\vspace{0.5cm}

\begin{Definition}
{\rm \cite{ka}} An almost Hermitian metric $\Omega$ on $(M^{2m}, J, \omega)$, where $m \geq 3$ is
said to be conformally balanced if it satisfies $d(\|\omega\|_{\Omega}\Omega^{m-1}) = 0$, where $\omega$ is a nowhere vanishing holomorphic $(m,0)$- form admitted by $M$, which means that $\omega$ satisfies $\bar{\partial}\omega = 0$. Then $(M, J, \omega, \Omega)$ is said to be an almost Calabi-Yau manifold with torsion. 
\end{Definition} 


\begin{Definition}
{\rm \cite{bric}} Let $M$ and $N$ be Riemannian manifolds, a smooth map $f: M \rightarrow N$ is harmonic if it is a critical point of the energy functional

$$
\mathbf{E}(f)=\frac{1}{2} \int_M\|\mathrm{~d} f\|^2 \mathrm{~d} v.
$$
Equivalently, $f$ solves the corresponding Euler-Lagrange equation $
\Delta f=0,$ where the Laplacian $\Delta f$ ( trace of the Hessain $\nabla^{2}f$), more commonly denoted by $\tau(f)$ and called the tension field, is a nonlinear operator generalizing the Riemannian Laplacian. We call ±holomorphic a map that is holomorphic or antiholomorphic.

\end{Definition}

\begin{Definition}
   {\rm \cite{bric}}  A smooth map $f: M \rightarrow N$, where $M$ is a complex manifold, is called pluriharmonic if for every $1$-dimensional complex submanifold $C \subseteq M$, the restriction $f\vert_C: C \rightarrow N$ is harmonic.\\
\noindent
Pluriharmonicity can alternatively be defined by :
$$
f \text { pluriharmonic } \quad \Leftrightarrow \quad \mathrm{d} \mathrm{d}^{\mathrm{c}} f=0 \text {. }
$$
\noindent
We also have the further characterizations:
$$
 f \text{ pluriharmonic } \quad \Leftrightarrow \quad \bar{\partial}\mathrm{d}^{1,0} f=0 \quad \Leftrightarrow \quad \partial\mathrm{d}^{0,1} f=0 \quad \Leftrightarrow \quad\left(\nabla^2 f\right)^{1,1}=0 .
$$
\end{Definition}
\begin{Remark}
 We have naturally written $\mathrm{d}=\mathrm{d}^{1,0}+\mathrm{d}^{0,1}=\partial+\bar{\partial}$. To make sense of the Hessian in the last one, choose any compatible Riemannian metric on $M$. Due to this characterization, pluriharmonic maps have sometimes been called $(1,1)$-geodesic maps.
\end{Remark}

\begin{Definition}
   {\rm \cite{allen}}  An $n$-dimensional vector bundle is a continuous function $p: E \rightarrow B$ together with a real vector space structure on $p^{-1}(b)$ for each $b \in B$, such that the following local triviality condition is satisfied: There is a cover of $B$ by open sets $U_\alpha$ for each of which there exists a homeomorphism $h_\alpha: p^{-1}\left(U_\alpha\right) \rightarrow U_\alpha \times \mathbb{R}^n$ taking $p^{-1}(b)$ to $\{b\} \times \mathbb{R}^n$ by a vector space isomorphism for each $b \in U_\alpha$. This $h_\alpha$ is called a local trivialization of the vector bundle. The space $B$ is called the base space (typically a manifold), $E$ is the total space which is a topological space, and the vector spaces $p^{-1}(b)$ are the fibers. 
\end{Definition}

\subsection{Some known results on astheno-K\"ahler manifold}
Here, we list the results based on the astheno-K\"ahler metric.

\begin{Lemma} {\rm\cite[Lemma 6]{jsty}} Let $(M,\Omega)$ be a compact astheno-K\"ahler manifold. Then, every holomorphic $1$-form on $M$ is closed.
\end{Lemma}

\noindent
The converse of the above lemma is not true. A counterexample  of manifold having no astheno-K\"ahler metrics with closed holomorphic $1$-form is given in \cite{at}  which is, let  $M$ be a compact quotient $ M = \Gamma \backslash G $ of a simply-connected Lie group
$G$ by a uniform discrete subgroup $\Gamma$, endowed with an invariant complex structure
$J$ and having no invariant astheno-K\"ahler metrics but any holomorphic $1$-form on $M$ is closed.


Jost and Yau gave the following results:
\begin{Lemma}{\rm\cite[Lemma 7]{jsty}} Let $N$ be a compact locally Hermitian symmetric space of noncompact type, and assume that the universal cover of $N$ does not have the upper half plane $\mathcal{H}$ as a global factor. Let $M$ be a compact astheno-K\"ahler manifold of dimension $m$. Let $f: M \rightarrow N$ be Hermitian harmonic. Then $f$ is pluriharmonic. If $f$ has real rank $2\dim_{\mathbb C} N$ at some point, then $f$ is $\pm$ holomorphic.
\end{Lemma}

\begin{Theorem}
{\rm\cite[Theorem 6]{jsty}}
Let $N$ be a compact locally Hermitian symmetric space of noncompact
type, without the upper half-plane $\cal H$ as a global factor of its universal cover, or let $N$ be a compact strongly negatively curved K\"ahler manifold. Let $M$ be a compact astheno-K\"ahler 
manifold. If $M$ is homotopy equivalent to $N$, then $M$ is $\pm$ biholomorphically equivalent to $N$.
\end{Theorem}

\begin{Remark}
The above theorem is the extension of Siu's rigidity theorem {\rm \cite[Theorem 2]{Siu}}.

\end{Remark}

\noindent
Fino and Tomassini also find the condition when astheno-K\"ahler becomes K\"ahler.
\begin{Theorem}
   {\rm \cite[Theorem 4.1]{at}} A conformally balanced astheno-K\"ahler structure $(J, g)$ on a compact manifold of complex dimension $m \geq 3$ whose Bismut connection has (restricted) holonomy contained in $SU(m)$ is necessarily K\"ahler and therefore it is a Calabi-Yau structure.
\end{Theorem}

\noindent
Latorre and Ugarte \cite{all} proved the following results:

\begin{Theorem}
    {\rm \cite[Corollary 2.3]{all}} Let $M$ be a homogeneous compact complex manifold of complex dimension $m \ge 3$ and let $g$ be an invariant Hermitian
metric on $M$. If $g$ is SKT or astheno-K\"ahler, then $g$ is $k$-th Gauduchon for any $1 \leq k \leq m - 1$.

\end{Theorem}

\begin{Theorem}
     {\rm \cite[Theorem 2.4]{all}} For each $m \ge 4$, there  is a non-K\"ahler compact complex manifold $M$ of complex dimension $m$ admitting a balanced
metric $\Tilde{g}$ and  an astheno-K\"ahler metric $g$, which is additionally $k$-th Gauduchon for any $1 \leq k \leq m - 1.$
\end{Theorem}

\noindent
In \cite{ak} Fino and Grantcherov proved that the twistor spaces of compact hyper-K\"ahler and negative quaternionic-K\"ahler manifolds do not admit astheno-K\"ahler metrics. They proved that compact semisimple Lie group of even dimension ($m > 6$) and rank
two endowed with its Samelson’s complex structure admits an astheno-K\"ahler non-SKT
metric. In contrast, its canonical SKT metric is not astheno-K\"ahler. In particular, the Lie groups $SU(3)$
and $G_2$ admit $SKT$ and astheno-K\"ahler metrics, which are different. 

\vskip .3cm 
\noindent 
Before listing further results, let us recall the following definitions from \cite{da}.
\begin{Definition}
 The Bott-Chern cohomology  of a complex manifold $X$ is the bi-graded algebra.
$$
H_{B C}^{\bullet, \bullet}(X):=\frac{\operatorname{ker} \partial \cap \operatorname{ker} \bar{\partial}}{\operatorname{im} \partial \bar{\partial}}.
$$
\noindent
Unlike in the case of the Dolbeault cohomology groups, for every $(p, q) \in \mathbb{N}^2$, the conjugation induces an isomorphism
$$
H_{B C}^{p, q}(X) \stackrel{\simeq}{\rightarrow} H_{B C}^{q, p}(X) .
$$

\end{Definition}

\begin{Definition}
 The Aeppli cohomology of a complex manifold $X$ is the bi-graded $H_{B C}^{\bullet, \bullet}(X)$-module
$$
H_A^{\bullet , \bullet}(X):=\frac{\operatorname{ker} \partial \bar{\partial}}{\operatorname{im} \partial+\operatorname{im} \bar{\partial}}   .
$$
\noindent
As for the Bott-Chern cohomology,  for every $(p, q) \in \mathbb{N}^2$, the conjugation induces the isomorphism
$$
H_A^{p, q}(X) \stackrel{\simeq}{\longrightarrow} H_A^{q, p}(X) .
$$
\end{Definition}

\noindent
 \textbf{Notation:} For every $(p, q) \in \mathbb{N}^2$, for every $k \in \mathbb{N}$, and for $\sharp \in$ $\{\bar{\partial}, \partial, B C, A\}$, we will denote
$$
h_{\sharp}^{p, q}:=\operatorname{dim}_{\mathbb{C}} H_{\sharp}^{p, q}(X) \in \mathbb{N} \quad \text { and } \quad h_{\sharp}^k:=\sum_{p+q=k} h_{\sharp}^{p, q} \in \mathbb{N}.
$$
Chiose and Rasdeaconu gave the following significant result on cohomologies on astheno-K\"ahler manifold:
\begin{Theorem}
  {\rm \cite[Theorem 1.1]{ro}}  On a compact astheno-K\"ahler manifold $(M,\Omega)$, the following inequalities
hold:
$$h_{BC}^{0,1}(M) \leq h_A^{0,1}(M) \leq h_{BC}^{0,1}(M) + 1.$$
\end{Theorem}

\begin{Theorem}
   {\rm \cite[Remark 4.2]{ro}} On an astheno-K\"ahler manifold $(M,\Omega)$, at least one of the
morphisms $$  H_{\Bar{\partial}}^{0,1}(M)\rightarrow H_{A}^{0,1}(M)  \hspace{.4cm} and \hspace{.4cm}  H_{BC}^{0,1}(M) \rightarrow H_{\Bar{\partial}}^{0,1}(M)$$is an isomorphism.
\end{Theorem}

\begin{Theorem}
    {\rm \cite[Remark 4.4]{ro}} On an astheno-K\"ahler manifold $(M,\Omega)$ of complex dimension $m$ one has $h_{BC}^{0,1} = h_A^{0,1}$ if and only if $\partial(e^{(m-1)f}\Omega^{m-1})$ is  $\partial\Bar{\partial}$- exact, where $f$ is a $C^\infty$ real function on $M$ such that $e^{(m-1)f}\Omega^{m-1}$ is $i\partial\Bar{\partial}-$ closed.
    \end{Theorem}

\begin{Theorem}
   {\rm \cite[Proposition 4.2]{ro}} Let $M$ be a compact complex astheno-K\"ahler manifold of dimension $m$, and let 
    $$ i_{p,q} : H_{BC}^{p,q}(M) \rightarrow H_{A}^{p,q}(M)$$
be the map induced by the identity. Then $i_{m-1,m-1}$ is surjective if and only if $i_{0,1} $ is injective.
\end{Theorem}

\begin{Theorem}
   {\rm \cite[Proposition 5.4]{ro}}  If $M_1$ and $M_2$ be two astheno-K\"ahler manifolds saturating the upper bound in $h_{BC}^{0,1}(M_1 \times M_2) \leq h_A^{0,1}(M_1 \times M_2) \leq h_{BC}^{0,1}(M_1 \times M_2) + 1$, then $M_1\times M_2$ does not admit astheno-K\"ahler metrics.
\end{Theorem}

\begin{Definition}
   {\rm \cite{vais}} A locally conformally K\"ahler (LCK)  manifold $(M,\Omega)$ is Vaisman if $\nabla\theta$ = 0, where $\nabla$ denotes the Levi-Civita connection associated to $\Omega$ and $\theta$ is the Lee form.
\end{Definition}

\begin{Theorem}
    {\rm \cite[Theorem 1.6]{ro}} Any compact Vaisman manifold of dimension at least three satisfies        $ h_A^{0,1}(M) = h_{BC}^{0,1}(M) + 1 $ and carries no astheno-K\"ahler metric and so being a Vaisman there exists no astheno-K\"ahler metrics on similarity Hopf manifolds of dimension at least three, diagonal Hopf manifold and standard Hopf manifold.
\end{Theorem}

\begin{Theorem}
   {\rm \cite[Corollary 6.4]{ro}} The class of astheno-K\"ahler manifolds is not invariant under modifications.
\end{Theorem}
In support of this, we can consider an example: let $M$ be the $3$-dimensional manifold constructed by Hironaka \cite{hiro}, which is a
proper modification of the projective space $\mathbb{CP}^3$ and which unlike $\mathbb{CP}^3$,  cannot carry an astheno-K\"ahler metric.

\vskip .3cm
\noindent
In  \cite{bis}, Biswas in his two-page paper gave a result that if $G/\Gamma$, where $\Gamma$ is a lattice in a connected complex Lie group $G$ such that $G/\Gamma$ is compact, admits an astheno-K\"{a}hler metric, then $G$ is abelian and $G/ \Gamma$ is a compact complex torus . The author used a lemma due to Jost and Yau \cite{jsty} that "any holomorphic $1$-form on a compact astheno-K\"{a}hler manifold is closed" as a key element in his proof.
\vskip .3cm
\noindent 
In \cite{fabio}, Podesta showed that for a suitable invariant complex structure
there exists an invariant Hermitian metric, which is Chern–Einstein and astheno-K\"ahler.

\vskip .2cm 
\noindent 
The main result of the paper reads as follows:
\begin{Theorem}
{\rm \cite[Theorem 3]{fabio}} Let $M=G / L$ be a manifold in the class $\mathrm{C}$ (special subclass of compact simply connected complex homogeneous manifolds that includes the
Calabi–Eckmann manifolds) endowed with an invariant complex structure $J_{a, b}(a, b \in \mathbb{R}, b \neq 0)$ and complex dimension $m$, where  $G$ is a compact connected Lie group and $L$ is some compact subgroup. Then $M$ is simply connected, non-K\"ahler and
there exist $a_o, b_o \in \mathbb{R}$ and an invariant metric, Hermitian concerning $J_{a_o, b_o}$, which is Chern-Einstein and astheno-K\"ahler.
\end{Theorem}

\noindent 
In \cite{ch}, Chen et al., in the second part of the paper, derived a characterization of $k$-Gauduchon metrics in terms of the Chern scalar curvature and a Riemannian type scalar curvature of $\Omega$ concerning the Chern connection, where in particular, a pluriclosed metric (i.e., $\partial\Bar{\partial}\Omega=0$) is $1$-Gauduchon. In contrast, an astheno-K\"ahler metric ( i.e., $\partial\Bar{\partial}\Omega^{m-2}=0$ ) is $(m-2)$-Gauduchon.
\vskip .2cm 
\noindent 
In \cite{ka}, Kawamura showed a relation between almost Gauduchon, astheno-K\"ahler and quasi-K\"ahler manifold. The author proved the following result:
\noindent
\begin{Theorem}
{\rm \cite[Theorem 1.1,1.2]{ka}} Suppose $(M^{2m},J,\omega,\Omega)$ is a compact almost-Calabi-Yau manifold with torsion for $m\geq 3$. Suppose  $\Omega$ is either almost-pluriclosed ( Gauduchon) or almost-astheno-K\"ahler. Then $\Omega$ is quasi-K\"ahler metric, where $\omega$ is a nowhere vanishing holomorphic $(m,0)$-form admitted by $M$. Here holomorphic means that $\omega$ satisfies $\bar{\partial}\omega = 0$.
\end{Theorem}

\noindent 
Calamai \cite{cs} provided a list of possible projectively flat metrics studied in \cite{jl} and showed that projectively flat astheno-K\"ahler metrics are globally conformally flat K\"ahler metrics.

\vskip .2cm 
\noindent 
 In \cite{st}, Sferruzza and Tomassini provided families of compact astheno-K\"ahler nilmanifolds endowed with left-invariant complex structure and studied the behaviour of blowup of compact complex manifolds endowed with astheno-K\"ahler metrics, satisfying certain extra differential conditions (i.e., $dd^c\Omega^{m-3}=0$ ) and showed that in this case the existence of astheno-K\"ahler metric is not preserved by blowup.
 \vskip .2cm
We will now present several findings related to nilmanifolds.
\begin{Remark}
  {\rm \cite{ag}}  Starting from late 1990, a "complex nilmanifold" means a quotient of a real nilpotent Lie group equipped with a left-invariant complex structure by the left action of a discrete, co-compact subgroup.  This definition is much more general; indeed, left-invariant complex structures are found on many even-dimensional nilpotent Lie groups that are not complex. The complex structure on a Kodaira surface is one such example.

    Complex structures on a nilmanifold have impeccable algebraic characterization. Let $G$ be a real nilpotent Lie group and $\mathfrak{g}$ be a Lie algebra. By Newlander-Nirenberg theorem \cite{new}, a complex structure on $G$ is the same as a sub-bundle $T^{1,0} G \subset T G \otimes_{\mathbb{R}} \mathbb{C}$ such that $\left[T^{1,0} G, T^{1,0} G\right] \subset T^{1,0} G$ and $T^{1,0} G \oplus \overline{T^{1,0} G}=T G \otimes_{\mathbb{R}}$ $\mathbb{C}$. The left-invariant sub-bundles in $T^{1,0} G$ are the same as subspaces $W \subset \mathfrak{g} \otimes_{\mathbb{R}} \mathbb{C}$, and the condition $\left[T^{1,0} G, T^{1,0} G\right] \subset T^{1,0} G$ is equivalent to $[W, W] \subset W$. Therefore, left-invariant complex structures on $G$ are the same as complex sub-algebras $\mathfrak{g}^{1,0} \subset \mathfrak{g} \otimes_{\mathbb{R}} \mathbb{C}$ satisfying $\mathfrak{g}^{1,0} \oplus \overline{\mathfrak{g}^{1,0}}=\mathfrak{g} \otimes_{\mathbb{R}} \mathbb{C}$.
\end{Remark}
\noindent Sferruzza \cite{sto} proved necessary cohomological conditions for the existence of curves of astheno-K\"ahler metrics
along curves of deformations starting from an initial compact complex manifold endowed with an
astheno-K\"ahler metric. Furthermore, he applied his results, providing obstructions to the existence
of curves of astheno-K\"ahler metrics on two different families of real 8-dimensional nilmanifolds
endowed with invariant nilpotent complex structures.
\begin{Theorem}
   {\rm \cite[Corollary 1.5]{ro}} Let $M$ be non-K\"ahler nilmanifold equipped with a nilpotent complex structure. If $M$ admits an astheno-K\"ahler metric, then $h_A^{0,1}(M) = h_{BC}^{0,1}(M) + 1.$
\end{Theorem}

\section{{Some known results on bundles on astheno-K\"ahler manifold}}
Let $f(x)=\sum\limits_{i=0}^{n}a_{i}x^{i}$ \ be a polynomial, where $a_{i}$
are nonnegative integers and $E$ be a holomorphic vector bundle on a compact complex manifold $X$. Let $f(E)=$ $\bigoplus\limits_{i=0}^{n}\left(
E^{\otimes i}\right) ^{\oplus a_{i}}$ be the vector bundle on $X$. An
holomorphic vector bundle $E$ on a compact complex manifold $X$ is called
finite \cite{bp} if there are two distinct polynomials $f,g$ with nonnegative integral
coefficients, such that the vector bundle $f(E)$ is isomorphic to $g(E)$. It
is known that a vector bundle $E$ is finite if and only if there is a finite
collection of vector bundles $\left\{ F_{j}\right\} _{j=1}^{m}$ and
nonnegative integers $\left\{ a_{i,j}\right\} _{j=1}^{m}$ such that

\[
E^{\otimes i}=\bigoplus\limits_{i=0}^{n}\left( F_{j}\right) ^{\oplus
a_{i,j}}
\]%
for all $i\geq 1$ \cite[p.35, Lemma 3.1]{Nori}.\\
\vspace{0.3cm}

\noindent Now, we give some results based on bundles.\\

\noindent
In \cite{bp}, Biswas and Pingali characterized finite vector bundles on a compact complex manifold admitting a Gauduchon and astheno-K\"ahler metric. The main result of the paper is as follows:
\begin{Theorem}
   {\rm \cite[Theorem 1.1]{bp}}  Let $M$ be a compact complex manifold that admits a Hermitian metric that is both Gauduchon and astheno-K\"ahler. Then a holomorphic vector bundle $E$ over $M$ is finite if and only if it corresponds to a representation of a finite quotient of the fundamental group of $M$, or equivalently, if and only if $E$ admits a flat connection, compatible with its holomorphic structure, that has a finite monodromy group.
    \end{Theorem}
  Before stating the next result, we need to state the following definition:

\begin{Definition}
{\rm \cite{cl} } Let $(X, \Omega)$ be a $m$-dimensional compact Hermitian manifold. A line bundle $L$ over $(X, \Omega)$ is said to be numerically effective if for every $\epsilon>0$, there exists a smooth metric $h_\epsilon$ on $L$ such that the curvature $\sqrt{-1} \Theta\left(L, h_\epsilon\right)=\sqrt{-1} \bar{\partial} \partial \log h_\epsilon \geq-\epsilon \Omega$.

\vskip .2cm 
\noindent
A Hermitian flat-line bundle is a numerically effective. A vector bundle, $E$ of rank $r \geq 2$, is said to be numerically effective flat if the anti-tautological line bundle $\mathcal{O}_E(1)$ on the projective bundle $P E$ is numerically effective flat. A vector bundle, $E$, is said to be numerically flat if both $E$ and its dual $E^*$ are numerically effective flat.
\end{Definition}

Chen \cite{cy} stated and proved the following result:
\begin{Theorem} Let $E$ be a pseudo-effective vector bundle on the Gauduchon and astheno-K\"ahler manifold $(M, \Omega)$ with vanishing first Chern number. Then, $E$ is numerically flat. 
\end{Theorem}

In \cite{cl}, Li et al. consider numerically flat vector bundles over compact non-K\"ahler manifolds, mainly astheno-K\"ahler manifolds by introducing the notion of approximate Hermitian flatness ( a holomorphic vector bundle $E$ over a compact complex manifold $X$ is approximate Hermitian flat if there exists a Hermitian metric $h_\epsilon$ such that the curvature satisfies $\sup _X\left|\Theta\left(L, h_\epsilon\right)\right|<\epsilon$, where $L$ is line bundle over a compact complex manifold $X$ and $\Theta\left(L, h_\epsilon\right)$ is the Chern form of an arbitrary Hermitian
 metric $h_\epsilon$ on $L$. ).\\

\vskip .2cm 
\begin{Theorem}
{\rm \cite{cl} } Let $(X, \widehat{\Omega})$ be an $m$-dimensional compact astheno-K\"ahler manifold, $\Omega$ a Gauduchon metric conformal to $\widehat{\Omega}$, and $E$ a holomorphic vector bundle over $X$.
\noindent
Then the following statements on $E$ are equivalent:\\
(1) $E$ is numerically flat.\\
(2) $E$ is numerically effective flat with $\operatorname{ch}_1(E) \cdot\left[\Omega^{m-1}\right]=0$.\\
(3) $E$ is $\Omega$-semistable with $\operatorname{ch}_1(E) \cdot\left[\Omega^{m-1}\right]=\operatorname{ch}_2(E) \cdot\left[\widehat{\Omega}^{m-2}\right]=0$.\\
(4) $E$ is approximate Hermitian flat.\\
(5) There exists a filtration

$$
0=E_0 \subset E_1 \subset E_2 \subset \cdots \subset E_l=E
$$
by subbundles whose quotients are Hermitian flat.
\end{Theorem}

\vskip .2cm 
\noindent 
 Biswas with Loftin  \cite{bl} contributed to the field of astheno-K\"ahler manifold. They proved the existence and uniqueness of a Harder-Narasimhan filtration for flat vector bundles over $M$, as well as a Bogomolov-type inequality for semistable flat vector bundles over $M$ under the assumption that the Gauduchon metric $g$ is astheno-K\"{a}hler.

\vskip .2cm 
\noindent 
Shen \cite{sz} studied the problem of Hitchin-Kobayashi correspondence for twisted holomorphic vector bundles over compact Gauduchon manifolds. The
main result proven is that a twisted holomorphic vector bundle over a compact Gauduchon manifold is semi-stable if and only if it admits an approximate Hermitian-Einstein
structure. The author derived a Bogomolov-type inequality for
a semi-stable twisted holomorphic vector bundle over a compact astheno-K\"ahler manifold.

\section{Examples}
Now, we give some examples of astheno-K\"ahler manifolds.

\begin{Example}
 {\rm \cite{km}}  Any product manifold of curves and surfaces is astheno-K\"ahler.
\end{Example}
\noindent
 Consider the product manifold $M=M_1 \times M_2$, where $M_i(i=1,2)$ be a $\left(2 m_i+1\right)$-dimensional compact normal almost contact metric manifold with the structure tensor fields $\left(\phi_i, \xi_i, \eta_i\right)$,  we consider an almost complex structure $J$ defined by
\begin{equation}\label{e1}
\begin{aligned}
J=\phi_1-\eta_2 \otimes \xi_1+\phi_2+\eta_1 \otimes \xi_2 
\end{aligned}
\end{equation}
 Here $M$ endowed with $J$ is a compact complex manifold of complex dimension $m=m_1+m_2+1$.
\noindent
Moreover, if $g_i$ is the compatible Riemannian metric on $M_i$ for each $i=1,2$, then the Riemannian product metric $g=g_1+g_2$ on $M$ is compatible with $J$, that is, $g$ is a Hermitian metric on $M$. 
\noindent
Then its K\"ahler form $\Omega$  is given by
\begin{equation}\label{e2}
\begin{aligned}
\Omega=\Phi_1+\Phi_2-2 \eta_1 \wedge \eta_2  
\end{aligned}
\end{equation}
where $\Phi_i$ denotes the fundamental 2-form on $M_i$ for each $i=1,2$.
\begin{Example}
    
\textbf{Product of two Sasakian manifolds}\\
{\rm \cite{km}} Let $\left(M_i, g_i\right)$ be a 3-dimensional compact Sasakian manifold for each $i=1,2$. Then the product manifold $M=M_1 \times M_2$ with the Hermitian structure  {\rm (\ref{e1})} and K\"ahler form {\rm (\ref{e2})} is astheno-K\"ahler.
\noindent
Since $M_1$ and $M_2$ are both Sasakian, we have
$$
d \Omega=-2\left(\Phi_1 \wedge \eta_2-\Phi_2 \wedge \eta_1\right)
$$
and
$$
d^{\mathrm{c}} \Omega=J d \Omega=-2\left(J \Phi_1 \wedge J \eta_2-J \Phi_2 \wedge J \eta_1\right)=2\left(\Phi_1 \wedge \eta_1+\Phi_2 \wedge \eta_2\right) .
$$
Using the fact that $\Omega$ and $\Phi_i$ are $J$-invariant, and $J \eta_1=\eta_2, J \eta_2=-\eta_1$, we get
$$
d d^{\mathrm{c}} \Omega=2\left(\Phi_1^2+\Phi_2^2\right)
$$
Since $\operatorname{dim}_{\mathbb{C}} M=m=3$, i.e., $m_i=1$ for each $i=1,2$, \ $\Phi_i^2=0$ on $M_i$, and hence
$$
d d^{\mathrm{c}} \Omega^{m-2}=d d^{\mathrm{c}} \Omega=2\left(\Phi_1^2+\Phi_2^2\right)=0 \quad \text { on } M \text {. }
$$
Therefore, the Hermitian structure  {\rm(\ref{e1})} on $M$ is astheno-K\"ahler.
\end{Example}

\noindent
{\bf  Note:}  If $\operatorname{dim}_{\mathbb{C}} M=m>3$ in above example, the Hermitian structure  (\ref{e1}) is not astheno-K\"ahler.
\begin{Example}

\textbf{ Product of  Sasakian manifold and Cosymplectic manifold}
\noindent\\
{\rm\cite{km}} Let $\left(M_1, g_1\right)$ be a 3-dimensional compact Sasakian manifold, and $\left(M_2, g_2\right)$ a compact cosymplectic manifold of dimension $\geq 3$. Then, the product manifold $M=M_1 \times M_2$ with the Hermitian structure  (\ref{e1})  is astheno-K\"ahler.\\
    \noindent
Since $M_1$ is Sasakian and $M_2$ is cosymplectic, we have
$$
d \Omega=-2 \Phi_1 \wedge \eta_2 \quad \text { and } \quad d^c \Omega=2 \Phi_1 \wedge \eta_1.
$$
By simple calculation, we have
$$
d d^{\mathrm{c}} \Omega=2 \Phi_1^2
$$
Since $m_1=1, \Phi_1^2=0$ on $M_1$, that is, $d d^{\mathrm{c}} \Omega=0$ and $d \Omega \wedge d^{\mathrm{c}} \Omega=0$ on $M$, and hence we obtain
$$
d d^{\mathrm{c}} \Omega \wedge \Omega+(m-3) d \Omega \wedge d^{\mathrm{c}} \Omega=0
$$
\noindent
Therefore, the Hermitian structure (\ref{e1}) on $M$ is astheno-K\"ahler.
\vskip .3cm
\noindent
{\bf Note:}  Let $\left(M_1, g_1\right)$ be a Sasakian manifold of dimension greater than $3$, and $\left(M_2, g_2\right)$ a cosymplectic manifold. Then the Hermitian structure (\ref{e1})is not astheno-K\"ahler.
\end{Example}
\begin{Example}

{\rm \cite{kmm}}   There exist astheno-K\"ahler structure (non-K\"ahler) on the Calabi-Eckmann manifold $M=S^{2 m_1+1} \times S^{2 m_2+1}$ with an almost complex structure $J$  given by Tsukada defined by
 \begin{equation}     
  J=\phi_1-\left(\frac{a}{b} \eta_1+\frac{a^2+b^2}{b} \eta_2\right) \otimes \xi_1+\phi_2+\left(\frac{1}{b} \eta_1+\frac{a}{b} \eta_2\right) \otimes \xi_2,
  \label{mori}\end{equation}
where $a, b \in \mathbb{R}$ and $b \neq 0$. The Hermitian metric $g$ on the complex manifold $(M, J)$ is given by
 $$g=g_1+g_2+a\left(\eta_1 \otimes \eta_2+\eta_2 \otimes \eta_1\right)+\left(a^2+b^2-1\right) \eta_2 \otimes \eta_2.$$
\noindent Then the K\"ahler form $\Omega$ on the Hermitian manifold $(M, J, g)$ is given by  $$\Omega=\Phi_1+\Phi_2-2 b \eta_1 \wedge \eta_2$$
where $\Phi_i$ denotes the fundamental $2$-form on $M_i$ for each $i=1,2$.

\end{Example}

\noindent {\bf Note:} 
When $a=0$ and $b=1$, the almost complex structure (\ref{mori}) coincides with A. Morimoto's complex structure (\ref{e1}).
\begin{Example}
   {\rm \cite[Theorem 1.3]{ro}} A cartesian product of two compact complex surfaces admits an
astheno-K\"ahler metric if and only if at least one of the surfaces admits a K\"ahler
metric.

\end{Example}

\begin{Example}
   {\rm \cite[Corollary 5.1]{ro}} The complex parallelizable Nakamura manifold does not admit astheno-K\"ahler metrics.
\end{Example}

\begin{Example}
   {\rm \cite[Corollary 5.3]{ro}} The Oeljeklaus-Toma manifolds of type $(s,t)$ with $s\ge 2$ and $s = 1$ do not admit astheno-K\"ahler metrics.
\end{Example}

\begin{Example}
{\bf By blow-ups and resolutions}\\
 \cite{at} Let $(M, J, g)$ be an astheno-K\"ahler manifold of complex dimension $m$ such that its fundamental $2$-form $\Omega$ satisfies $$\partial\Bar{\partial}\Omega = 0 ,\partial\Bar{\partial}\Omega^2 = 0. $$ Then both the blow-up
$ \Tilde{M}_p $
of $M$ at a point $p \in M$ and the blow-up $ \Tilde{M}_Y $ of $M$ along a compact complex
submanifold $Y$ admit an astheno-K\"ahler metric satisfying $\partial\Bar{\partial}\Omega = 0 ,\partial\Bar{\partial}\Omega^2 = 0.$ 
Thus, it is possible to construct new examples of astheno K\"ahler manifolds by blowing up a given astheno-K\"ahler manifold M (satisfying $\partial\Bar{\partial}\Omega = 0,\ \partial\Bar{\partial}\Omega^2 = 0 $) at one or more points or along a compact complex submanifold.
\end{Example}
\begin{Example}{\bf  By twist construction }\\
 \cite{at} Applying the twist construction {\rm \cite[Prop. 4.5]{as}} to $8$-dimensional astheno-K\"ahler manifolds with torus action, one can get new simply connected astheno-K\"ahler manifolds.
Let $\left(N^6, J,\right)$ be a $6$-dimensional simply connected compact complex manifold with a Hermitian structure $g$, which is strong KT and standard. Consider the product $M^8=N^6 \times \mathbb{T}^2$, where $\mathbb{T}^2$ is a $2$ -torus with an invariant K\"ahler structure. Then $M^8$ is astheno-K\"ahler and strong KT with torsion $c$ supported on $N^6$.

\end{Example}
\begin{Example}  {\bf $8$-dimensional nilmanifolds}\\ {\rm\cite[Example 2.3]{at}} 
Let $\left\{\eta^1, \ldots, \eta^4\right\}$ be the set of complex forms of type $(1,0)$ such that 
$$
\left\{\begin{aligned}
d \eta^j= & 0, \quad j=1,2,3 \\
d \eta^4= & a_1 \eta^1 \wedge \eta^2+a_2 \eta^1 \wedge \eta^3+a_3 \eta^1 \wedge \bar{\eta}^1+a_4 \eta^1 \wedge \bar{\eta}^2+a_5 \eta^1 \wedge \bar{\eta}^3 \\
& +a_6 \eta^2 \wedge \eta^3+a_7 \eta^2 \wedge \bar{\eta}^1+a_8 \eta^2 \wedge \bar{\eta}^2+a_9 \eta^2 \wedge \bar{\eta}^3+a_{10} \eta^3 \wedge \bar{\eta}^1 \\& +a_{11} \eta^3 \wedge \bar{\eta}^2+a_{12} \eta^3 \wedge \bar{\eta}^3
\end{aligned}\right.
$$
span the dual of a 2-step nilpotent Lie algebra $\mathfrak{n}$, depending on the complex parameters $a_1, \ldots, a_{12}$ and define an integrable almost complex structure $J$ on $\mathfrak{n}$.
 Let $a_1, \ldots, a_{12} \in \mathbb{Q}[i]$ and $(M=\Gamma \backslash N, J)$ be the corresponding compact complex nilmanifold of real dimension 8, where $N$ is connected Lie group with Lie algebra $\mathfrak{n}$ and $\Gamma$ is uniform discrete subgroup, then the Hermitian metric
$$
g=\frac{1}{2} \sum_{j=1}^4 \eta^j \otimes \bar{\eta}^j+\bar{\eta}^j \otimes \eta^j
$$
is astheno-K\"ahler if and only if
$$
\begin{gathered}
\left|a_1\right|^2+\left|a_2\right|^2+\left|a_4\right|^2+\left|a_5\right|^2+\left|a_6\right|^2+\left|a_7\right|^2+\left|a_9\right|^2+ \\
\left|a_{10}\right|^2+\left|a_{11}\right|^2=2 \Re \mathfrak{e}\left(a_3 \bar{a}_8+a_3 \bar{a}_{12}+a_8 \bar{a}_{12}\right) .
\end{gathered}
$$
If, in addition $a_8=0$ and $\left|a_4\right|^2+\left|a_{11}\right|^2 \neq 0$, then the astheno-K\"ahler metric $g$ is not SKT. Moreover, if $a_8=0$, the astheno-K\"ahler metric $g$ is SKT if and only if $a_1=a_4=a_6=a_7=a_9=a_{11}=0$.

\end{Example}
In their study, Latorre and Ugarte \cite{all} presented instances of non-K\"ahler compact complex manifolds supporting balanced and astheno-K\"ahler metrics. They provided constructions of compact complex non-K\"ahler manifolds with complex dimension $m \geq 4$, which possess a balanced metric and an astheno-K\"ahler metric that is $k$-Gauduchon for any $k$ within the range $\{1,\ldots, m-1\}$. Additionally, they offered examples of nilmanifolds characterized by a Lie algebra that is a product of a $(2n + 1)$-dimensional Heisenberg Lie algebra and a real line featuring an invariant abelian complex structure. Furthermore, they also explored examples on nilmanifolds with an invariant non-abelian complex structure specifically in complex dimension 4.
\begin{Example}
     {\rm \cite[Example 4.2]{ak}} The total space of a principal $T^2$-bundle over a 6-dimensional torus  admits a balanced and astheno-K\"ahler metric. Let $\pi: M \rightarrow T^6$ be the principal $T^2$ bundle over $T^6$ with characteristic classes
$$
a_1=d z_1 \wedge d \bar{z}_1+d z_2 \wedge d \bar{z}_2-2 d z_3 \wedge d \bar{z}_3, \quad a_2=d z_2 \wedge d \bar{z}_2-d z_3 \wedge d \bar{z}_3,
$$
where $\left(z_1, z_2, z_3\right)$ are complex coordinates on $T^6$. Consider on $T^6$ the standard complex structure and let
$$
F_1=d z_1 \wedge d \bar{z}_1+d z_2 \wedge d \bar{z}_2+d z_3 \wedge d \bar{z}_3, \quad F_2=d z_1 \wedge d \bar{z}_1+d z_2 \wedge d \bar{z}_2+5 d z_3 \wedge d \bar{z}_3,
$$
then the $a_i$ 's are traceless with respect to $F_1$ and $\left(a_1^2+a_2^2\right) \wedge F_2=0$. Let $\theta^j$ be connection 1-forms such that $d \theta^j=\pi^* a_j$ and define
$$
\Omega_1=\pi^* F_1+\theta^1 \wedge \theta^2, \quad \Omega_2=\pi^* F_2+\theta^1 \wedge \theta^2 .
$$
The $2$-forms $\Omega_1$ and $\Omega_2$ define a balanced metric $g_1$ and an astheno K\"ahler metric $g_2$ on $M$, respectively, compatible with the integrable complex structure so that the projection map $\pi$ is holomorphic. $M$ can be alternatively described as the $2$-step nilmanifold $G / \Gamma$, where $G$ is the $2$-step nilpotent Lie group with structure equations
$$
\left\{\begin{array}{l}
d e^j=0, \quad j=1, \ldots, 6 \\
d e^7=e^1 \wedge e^2+e^3 \wedge e^4-2 e^5 \wedge e^6, \\
d e^8=e^3 \wedge e^4-e^5 \wedge e^6
\end{array}\right.
$$
and $\Gamma$ is a co-compact discrete subgroup, endowed with the invariant complex structure $I$ such that $ Ie_1 = e_2, Ie_3 = e_4, Ie_5 = e_6, Ie_7 = e_8 $. In this setting the 2-form $ e^1\wedge e^2 + e^3 \wedge e^4 + e^5 \wedge e^6 + e^7 \wedge e^8 $ defines a balanced metric and $ e^1\wedge e^2 + e^3 \wedge e^4 + 5e^5 \wedge e^6 + e^7 \wedge e^8 $ gives an astheno-K\"ahler metric.

\end{Example}

\begin{Example}
 {\rm \cite[Proposition 5.7]{ak}} The homogeneous space $SU(5)/T^2$ ( flag manifold) for appropriate action of $T^2$ is
simply connected. It has an invariant complex structure which admits both balanced
and astheno-K\"ahler metrics but does not admit any $SKT$ metric.

\end{Example}

\begin{Example} In {\rm\cite{kmm}} Matsuo showed the existence of astheno-K\"ahler metrics on Calabi-Eckmann
manifolds. Since Calabi-Eckmann manifolds are principal $T^2$-bundles over $\mathbb{CP}^n\times\mathbb{CP}^m$, the construction of astheno-K\"ahler structures on torus bundles over K\"ahler manifolds. 
\end{Example}
In \cite{ak}, Fino et al. generalised Matsuo’s result to principal torus fibrations over
 compact K\"ahler manifolds. Using this idea Fino et al. gave a new example, which is 
 \begin{Example}
 Let $ \pi: P \rightarrow M$ be an $n$-dimensional principal torus bundle over a
 K\"ahler manifold $(M,J,\Omega)$ equipped with $2$ connections $1$-forms $\theta_{1},\theta_{2}$ whose curvatures
$ Q_{1}, Q_{2}$ are of type $(1,1)$ and are pull-backs from forms $\alpha_{1}$ and $\alpha_{2}$ on $M$. Let $I$ be the
 complex structure on $P$ defined as the pull-back of $J$ to the horizontal subspaces and
 as $I(\theta_{1}) = \theta_{2}$ along vertical directions. Let $\Phi = \pi^{\ast}(\Omega) + \theta_{1} \wedge \theta_{2}$; then
 $dd^{c}\Phi^{k} = k(Q_{1}^{2}+Q_{2}^{2})\wedge (\pi^{\ast}(\Omega^{k-1})), $
  where $1 \leq k \leq n-2$.
 In particular, if $\alpha_{1}^{2}+\alpha_{2}^{2}\wedge \Omega^{n-3} = 0$, then $\Phi$ is astheno-K\"ahler and if a torus bundle with $2$-dimensional fiber over a K\"ahler base admits an SKT metric, then it is astheno-K\"ahler.
 \end{Example}


\begin{Example}
    {\rm \cite[Proposition 5.6]{ro}} Let $S$ be a compact complex surface, $M$ a compact K\"ahler manifold. If $N$ denotes the blow-up of $M \times S$ along a smooth submanifold, then $N$ is an astheno-K\"ahler manifold.
\end{Example}


\section{Equations and  estimates of their solutions on astheno-K\"ahler manifold}
This section explores the equation on compact Hermitian manifolds and demonstrates the existence of solutions for the equation on astheno-K\"ahler manifolds.
\subsection{The Fu-Yau equation on compact astheno-K\"ahler manifolds}

Let $(M, \Omega)$ be an $m$-dimensional compact K\"ahler manifold. As a reduced generalized Strominger system in higher dimensions, Fu and Yau introduced the following fully nonlinear equation for $\varphi $, which is usually called Fu-Yau equation \cite{dhp, mg}

\begin{eqnarray}
\sqrt{-1} \partial \bar{\partial}\left(e^{\varphi} \Omega-\alpha e^{-\varphi} \rho\right) \wedge \Omega^{m-2} 
+m \alpha \sqrt{-1} \partial \bar{\partial} \varphi \wedge \sqrt{-1} \partial \bar{\partial} \varphi \wedge \Omega^{m-2}+\mu \frac{\Omega^m}{m !}=0,
\label{eq-aaaa}
\end{eqnarray}
where $\alpha$ is a non-zero constant called the slope parameter, $\rho$ is a real smooth $(1,1)$ form, $\mu$ is a smooth function. For $\varphi$, we can impose the elliptic condition
$$
\tilde{\Omega}=e^{\varphi} \Omega+\alpha e^{-\varphi} \rho+2 m \alpha \sqrt{-1} \partial \bar{\partial} \varphi \in \Gamma_2(M)
$$
and the normalization condition
\begin{equation}
\left\|e^{-\varphi}\right\|_{L^1}=A
\label{eq-aa}
\end{equation}
where
$$
\Gamma_2(M)=\left\{\alpha \in A^{1,1}(M) \left|  \ \frac{\alpha^1 \wedge \Omega^{m-1}}{\Omega^m}>0, \frac{\alpha^2 \wedge \Omega^{m-2}}{\Omega^m}>0\right\}\right.
$$
and $A^{1,1}(M)$ is the space of smooth real $(1,1)$ forms on $M$.
When $m=2$, (\ref{eq-aaaa}) is equivalent to the Strominger system on a toric fibration over a $K 3$ surface constructed by Goldstein and Prokushki \cite{Gold}, which was solved by Fu and Yau for $\alpha>0$ and $\alpha<0$ in \cite{jxm,jxst}.
\noindent
In case of $\alpha<0$, Phong et al. \cite{dhp}  proved the existence of solutions of (\ref{eq-aaaa}) such that the condition (\ref{eq-aa}) is replaced by
$$
\left\|e^{\varphi}\right\|_{L^1}=\frac{1}{A} \gg 1.
$$




Chu et al. \cite{jlx} proved the following result:
\begin{Theorem} 
 Let $(M, \Omega)$ be an $m$-dimensional compact astheno-K\"ahler manifold. Then there exist constants $A_0, C_0, \delta_0, M_0$ and $D_0$ depending only on $\alpha, \rho, \mu$ and $(M, \Omega)$ such that for any $A \leqslant A_0$, there exists a unique solution $\varphi$ of the Fu-Yau equation and
satisfying the elliptic condition and the normal condition $$
e^{-\varphi} \leqslant \delta_0,\ \ \ |\partial \bar{\partial} \varphi|_g \leqslant D,\ \ \  D_0 \leqslant D \ \ \ \text { and } A \leqslant \frac{1}{C_0 M_0 D}.
$$
\end{Theorem}
{\bf Remark:}

\begin{itemize}
    \item They used the Moser iteration to do $\mathrm{C}^0$-estimate for solutions $\varphi$ of the Fu-Yau equation and also gave a sketch of $\mathrm{}{C}^1$,$\mathrm{C}^2$ estimates of $\varphi$.

    \item Furthermore, they improved a result for the case $\alpha \leq 0$ and provided another proof to derive a prior $\mathrm{C}^1,\mathrm{C}^2$ estimate for $\varphi$ of the Fu-Yau equation in the case of $\alpha \leq 0$, but without the restriction condition.  

\end{itemize}
Kawamura  \cite{kaw} investigated the Fu-Yau equation on compact almost astheno-K\"ahler manifolds and showed an a
priori $C^0$-estimate for a smooth solution of the equation.
\subsection{Monge-Amp\`ere Equation on compact Astheno-K\"ahler manifolds}
The complex Monge-Amp\`ere equation
$$
(\Omega+\sqrt{-1} \partial \bar{\partial} u)^m=e^F \Omega^m, \quad \Omega+\sqrt{-1} \partial \bar{\partial} u>0,
$$
on a compact K\"ahler manifold $(M, \Omega)$ of dimension $m$ was solved by Yau in the 1978 \cite{yo} and has consistently been influenced by its pervasive role in K\"ahler geometry since its inception. Here $F$ is a given smooth function on $M$, normalized so that $\int_M e^F \Omega^m=\int_M \Omega^m$. Yau's Theorem, conjectured in the 1950s by Calabi, is that the Monge-Amp\`ere equation has a unique solution $u$ with $\sup _M u=0$.
\vskip.5cm
\noindent
For a general Hermitian metric $\Omega$, the complex Monge-Amp\`ere equation was solved in full generality by the authors in \cite{tv}. In this case, there exists a unique pair $(u, b)$ with $u$ a smooth function satisfying $\sup _M u=0$ and $b$ a constant such that
$$
(\Omega+\sqrt{-1} \partial \bar{\partial} u)^m=e^{F+b} \Omega^m, \quad \Omega+\sqrt{-1} \partial \bar{\partial} u>0 .
$$
This equation applies to studying cohomology classes and notions of positivity on complex manifolds.

Tosatti and Weincove \cite{vt} showed unique smooth solutions to the Monge-Amp\`ere equation for $(m-1)$-plurisubharmonic functions on Hermitian manifolds. They  obtained Calabi-Yau theorems for Gauduchon and strongly Gauduchon
metrics on a class of non-K\"ahler manifolds those satisfying the astheno-K\"ahler condition.
They discussed another Monge-Amp\`ere equation  (introduced by Popovici \cite{popo}) and
showed that the full Gauduchon conjecture can be reduced to a second
order estimate of Hou-Ma-Wu type.

\vspace{0.1cm}

\noindent The main result of the authors is given below: 
\begin{Theorem}
Let $M$ be a compact complex manifold equipped with an astheno-K\"ahler metric $\Omega$. Let $\Omega_0$ be a Gauduchon (resp. strongly Gauduchon) metric on $M$ and $F^{\prime}$ a smooth function on $M$. Then there exists a unique constant $b^{\prime}$ and a unique Gauduchon (resp. strongly Gauduchon) metric, which we write as $\Omega_u$, with
$$
\Omega_u^{m-1}=\Omega_0^{m-1}+\sqrt{-1} \partial \bar{\partial} u \wedge \Omega^{m-2}
$$
for some smooth function $u$.
    
\end{Theorem} 




\subsection{Hermitian-Yang-Mills flow}

Let $\left(E, \bar{\partial}_E\right)$ be a rank $r$ holomorphic vector bundle over a $m$-dimensional compact complex manifold $M$ with Hermitian metric $g$ and associated $(1,1)$-form $\Omega$. Let $H_0$ be the initial Hermitian metric on $E$ and $\mathcal{A}_{H_0}$ be the space of connections of $E$ compatible with ${H_0}$.  Let $F_A$ denote the curvature of the connection $A\in\mathcal{A}_{H_0}$ and $\Lambda_\Omega$ denote the adjoint of the operation $\eta \mapsto \eta \wedge \Omega$. 

In  \cite{nz}, Nie and Zhang studied the limiting behaviour of solutions of the Hermitian-Yang-Mills flow on $\mathcal{A}_{H_0}$,
\begin{equation}\label{e22}
\begin{aligned}
\frac{\partial A}{\partial t} & =i\left(\bar{\partial}_A-\partial_A\right) \Lambda_\Omega F_A, \\
A(0) & =A_0.
\end{aligned}
\end{equation}

\noindent
When $(M, g)$ is K\"ahler, the flow equations \eqref{e22} reduce to the Yang-Mills flow,
\begin{equation}\label{e3}
\begin{aligned}
\frac{\partial A}{\partial t}=-d_A^* F_A .
\end{aligned}
\end{equation}
\vskip .3cm
\noindent
In \cite{nz}, the authors proved that 
\begin{Theorem}
    
Let $(M, g)$ satisfy the Gauduchon condition $\partial \bar{\partial} \Omega^{m-1}=0$ and the astheno-K\"ahler condition $\partial \bar{\partial} \Omega^{m-2}=0$  and let $A(t)$ be a solution of \eqref{e22}. Then for every sequence $\left\{t_k\right\}_{k=1}^{\infty} \subset \mathbb{R}$ with $\underset{k \rightarrow \infty}{{\rm lim}} t_k=\infty$, there is a subsequence $t_{k_j}$ and a closed set $\Sigma \subset M$ of Hausdorff codimension at least four such that $A(t_{k_j})$ converges modulo gauge transformations to a connection $A_{\infty}$ on a Hermitian vector bundle $\left(E_{\infty}, H_{\infty}\right)$ in $C_{\text {loc }}^{\infty}$ topology on $M \backslash \Sigma$ , where $A_{\infty}$ satisfies
$$
d_{A_{\infty}} \Lambda_\Omega F_{A_{\infty}}=0.
$$
\end{Theorem}
\noindent{\bf Acknowledgement:} We would like to express our sincere appreciation to anonymous reviewers for their thoughtful and constructive feedback, which greatly contributed to enhancing the quality of the paper.


\begin{thebibliography}{99}


 \bibitem{scyau} W. Ballmann, Lectures on K\"{a}hler manifolds, ESI Lectures in mathematics and physics, (2006).


 \bibitem{bl} I. Biswas and J. Loftin, Hermitian-Einstein connections on principal bundles over flat affine manifolds, Int. Jour. Math., 23 (2012), no. 4, 1250039, 23 pp.

\bibitem{bis} I. Biswas, On Hermitian structure on $G/\Gamma$, Bull. Sci. Math., 137 (2013), 716–717.

\bibitem{bp} I.Biswas and V. P. Pingali, A characterization of finite vector bundles on Gauduchon astheno-K\"ahler manifolds, Épijournal de Géométrie Algébrique, 2 (2018), Art. 6, 13pp.

\bibitem{da} R. Bott and S.S. Chern, Hermitian vector bundles and the equidistribution of the zeroes of their holomorphic sections, Acta Math., 114 (1965), no. 1, 71-112.


 \bibitem{cs} S. Calamai, Positive projectively flat manifolds are locally conformally flat-K\"ahler Hopf manifolds, Pure Appl. Math. Q., 17 (2021), no.3, 1139-1154.

\bibitem{ch} H. Chen,  L. Chen and X. Nie, Chern-Ricci curvatures, holomorphic sectional curvature and Hermitian metrics, Sci. China Math., 64 (2021), no.4, 763-780.

\bibitem{cy} Y. Chen, A note on pseudo-effective vector bundles with vanishing first Chern number over non-K\"ahler manifold, C. R. Math. Acad. Sci. Paris, 359 (2021), 523-531.


\bibitem{ro} I. Chiose and R. Rasdeaconu, Remarks on astheno-K\"ahler manifolds, bott-chern and aeppli cohomology groups, Ann. Glob. Anal. Geom., 63 (2023), no. 24.



\bibitem{jlx} J. Chu , L. Huang and X. Zhu, The Fu–Yau equation on compact astheno-K\"ahler manifolds, Adv. Math., 346 (2019), 908-945.

\bibitem{ag} A. Fino, G. Grantcharov and M. Verbitsky, Algebraic dimension of complex nilmanifolds, J. Math. Pures Appl., 118 (2018), 204–218.



\bibitem{ak} A. Fino, G. Grantcharov and L. Vezzoni, Astheno–K\"ahler and Balanced Structures on Fibrations, Int. Math. Research Notices, 2019  (2019), no. 22, 7093–7117.
\bibitem{at}  A. Fino and A. Tomassini, On astheno-K\"ahler metrics, Jour. Lond. Math. Soc., 83 (2011), no. 2, 290-308.

\bibitem{jxm} J.-X. Fu and S.-T. Yau, A Monge-Amp\`ere-type equation motivated by string theory, Comm. Anal. Geom., 15 (2007), no. 1, 29–75.

\bibitem{jxst} J.-X. Fu and S.-T. Yau, The theory of superstring with flux on non-K\"ahler manifolds and the complex Monge-Amp\`ere equation, Jour. Diff. Geom., 78 (2008), no. 3, 369–428.

\bibitem{fu} J. Fu, Z. Wang and D. Wu, Semilinear equations, the function, and generalized Gauduchon metrics, Jour. Eur. Math. Soc., 15 (2013), 659-680.

\bibitem{mg} M. Garcia-Fernandez, Lectures on the Strominger system, Travaux math\'ematiques, 24 (2016), 7-61.

\bibitem{gaud} P. Gauduchon, La $1$-forme de torsion d’une vari\'et\'e hermitienne compacte, Math. Ann. 267 (1984), no. 8, 495–518.

\bibitem{Gold} E. Goldstein and S. Prokushkin, Geometric model for complex non-K\"ahler manifolds
with $SU(3)$ structure, Comm. Math. Phys., 251 (2004), no. 1, 65–78.

\bibitem{allen} A. Hatcher, Vector Bundles and K-theory, Version 2.1, (May 2009).

\bibitem{hiro} H. Hironaka, An example of a non-K\"ahlerian complex-analytic deformation of K\"ahlerian complex structures, Ann. Math., 75 (1962), no. 2, 190–208.

\bibitem{jsty} J. Jost and S.-T. Yau, A nonlinear elliptic system for maps from Hermitian to Riemannian manifolds and rigidity theorems in Hermitian geometry, Acta Math., 170 (1993), 221–254.

\bibitem{jst} J. Jost and S.T. Yau, Correction to A nonlinear elliptic system for maps 
from Hermitian to Riemannian manifolds and rigidity theorems in Hermitian geometry, Acta Math., 173 (1994), 307.


\bibitem{ka} M. Kawamura,  On the conformally balanced condition on almost Hermitian manifolds and the quasi-K\"ahlerity, J. Geom., 112 (2021), no.2,  paper no. 20, 13 pp.


\bibitem{kaw} M. Kawamura, An a priori $C^0$-estimate for the Fu-Yau
equation on compact almost astheno-K\"ahler manifolds, Complex Manifolds, 9 (2022), no. 1, 223-237.

\bibitem{kobayashi} S. Kobayashi, Kan\^o Memorial Lectures 5, Differential geometry of Complex vector bundles, Iwanami Shoten, Publishers
and Princeton University Press (1987).




\bibitem{all} A. Latorre and L. Ugarte, On non-K\"ahler compact complex manifolds with balanced and astheno-K\"ahler metrics, C. R. Math. Acad. Sci. Paris, 355 (2017), no. 1, 90–93.  

\bibitem{cl} C. Li, Y. Nie and X. Zhang, Numerically flat holomorphic bundles over non-K\"ahler manifolds, Jour. Reine Angew. Math., 790 (2022), 267–285.

\bibitem{nijen} A. Nijenhuis, $X_{n-1}$ forming sets of eigenvectors, Nederl. Akad. Wet., Proc., Ser. A, 54 (1951), 200-212; Indagationes Math. 13 (1951), 200-212.
 
\bibitem{jl} J. Li, S.T. Yau and F. Zheng, On projectively flat Hermitian manifolds, Comm. Anal. Geom., 2 (1994), no. 1, 103–109.  


\bibitem{bric} B. Loustau, Harmonic maps from K\"ahler manifolds, https://arxiv.org/abs/2010.03545v1, (2020).


\bibitem{km} K. Matsuo and T. Takahashi, On compact astheno-K\"ahler manifolds, Colloq. Math., 89 (2001), 213– 221.

\bibitem{kmm} K. Matsuo, Astheno-K\"{a}hler structures on Calabi–Eckmann manifolds, Colloq. Math., 115 (2009), 33–39.


\bibitem{new} A. Newlander and L. Nirenberg, Complex analytic coordinates in almost complex manifolds, Ann. Math., 65 (1957), 391-404.

\bibitem{nz} Y. Nie and X. Zhang, The limiting behaviour of the Hermitian-Yang-Mills
flow over compact non-K\"ahler manifolds, Sci. China Math., 63 (2020), no. 7, 1369-1390. 

 \bibitem{Nori} M. V. Nori, On the representations of the fundamental group, Compositio Math., 33 (1976), no. 1, 29–41.
 
\bibitem{dhp}D.H. Phong, S. Picard and X. Zhang, The Fu-Yau equation with negative slope
parameter, Invent. Math., 209 (2017), no. 2, 541–576.

\bibitem{fabio} F. Podest\`a, Homogeneous Hermitian manifolds and special metrics, Transformation Groups, 23 (2018), no. 4, 1129–1147.

 \bibitem{popo} D. Popovici, Aeppli cohomology classes associated with Gauduchon metrics on compact complex manifolds, Bull. Soc. Math. Fr., 143 (2015), no. 4, 763-800.
 

\bibitem{sz} Z.Shen, Semi-stable twisted holomorphic vector bundles over Gauduchon manifolds, Bull. Sci. Math., 187 (2023), paper No. 103288, 21 pp.


\bibitem{Siu} Y.T. Siu, The complex analyticity of harmonic maps and the strong rigidity of compact K\"ahler manifolds, Annals Math., 112 (1980), 73-111.

\bibitem{Streets} J. Streets, and G. Tian, A parabolic flow of pluriclosed metrics, Int. Math. Res. Not. IMRN, 2010 (2010), no.16, 3101–3133.

\bibitem{as} A. Swann, Twisting Hermitian and hypercomplex geometries, Duke Math. Jour., 155 (2010), no. 2, 403-431.


\bibitem{sto} T. Sferruzza, Deformations of astheno-K\"ahler metrics, Complex Manifolds, 10 (2023), no. 1, 19 pp.

\bibitem{st} T. Sferruzza and A. Tomassini, On cohomological and formal properties of strong K\"ahler
with torsion and astheno-K\"ahler metrics, Math. Z., 304 (2023), paper no. 55, 27pp.


   

\bibitem{tv} V. Tosatti and B. Weinkove, The complex Monge-Amp\`ere equation on compact Hermitian manifolds, Jour. Amer. Math. Soc., 23 (2010), no.4, 1187–1195.



\bibitem{vt} V. Tosatti and B. Weinkove, Hermitian metrics, $(n-1,n-1)$ form and Monge-Amp\`ere Equations, Journal für die reine und angewandte Mathematik (Crelles Journal), 2019 (2019), no. 755, 67-101.

\bibitem{vais} I. Vaisman, On locally conformal almost Kaehler manifolds, Isr. J. Math., 24 (1976), 338-351.
\bibitem{Yano} K. Yano and M. Kon, Structures on manifolds, Series in Pure Mathematics, World Scientific, Singapore, Distr. by John Wiley and Sons Ltd., Chichester. IX, 3 (1984), 508 pp.

\bibitem{yo} S.T. Yau, On the Ricci curvature of a compact K\"ahler manifold and the complex Monge-Amp\`ere equation, I, Comm. Pure Appl. Math., 31 (1978), no.3, 339–411.
\end{thebibliography}
\end{document}